\documentclass[12pt]{article}
\usepackage{amsmath, amssymb, amsthm, amsfonts, enumerate}
\usepackage{graphicx}

\newtheorem{theorem}{Theorem}
\newtheorem{lemma}[theorem]{Lemma}
\newtheorem{proposition}[theorem]{Proposition}
\newtheorem{remark}[theorem]{Remark}
\newtheorem{cor}[theorem]{Corollary}

\def\rr{{\mathbb R}}

\def\sik{{\rr}^2}

\def\cc{{\mathbb C}}

\def\qq{{\mathbb Q}}

\def\se{\setminus}

\def\al{\alpha}
\def\be{\beta}
\def\ga{\gamma}

\def\de{\delta}
\def\De{\Delta}
\def\ep{\varepsilon}
\def\si{\sigma}

\def\ph{\phi}
\def\Ph{\Phi}

\def\cd{\cdot}
\def\stb{,\ldots ,}
\def\ha{\frac{1}{2}}

\def\del{\partial}

\def\msk{\medskip}
\def\bsk{\bigskip}
\def\noi{\noindent}

\def\({\left(}
\def\){\right)}
\def\[{\left[}
\def\]{\right]}

\def\ol{\overline}

\def\sumi0n{\sum_{i=0}^n}

\def\proof{\noi {\bf Proof.} }

\def\gcd{{\rm gcd} \, }

\begin{document}

\title{Irregular tilings of regular polygons with similar triangles}

\author{M. Laczkovich}

\footnotetext[2]{{\bf MR subject classification:} 52C20}
\footnotetext[3]{The author was supported by the Hungarian National
Foundation for Scientific Research, Grant No. K124749.}

\maketitle

\begin{abstract}
We say that a triangle $T$ tiles a polygon $A$, if $A$ can be dissected
into finitely many nonoverlapping triangles similar to $T$. We show that if
$N>42$, then there are at most three nonsimilar triangles $T$ such that the
angles of $T$ are rational multiples of $\pi$ and $T$ tiles the regular $N$-gon.

A tiling into similar triangles is called regular, if the pieces have
two angles, $\al$ and $\be$, such that at each
vertex of the tiling the number of angles $\al$ is the same as that of $\be$.
Otherwise the tiling is irregular. It is known that for every regular polygon
$A$ there are infinitely many triangles that tile $A$ regularly. We show that
if $N>10$, then a triangle $T$ tiles the regular $N$-gon
irregularly only if the angles of $T$ are rational multiples of $\pi$.
Therefore, the numbers of triangles tiling the regular $N$-gon irregularly
is at most three for every $N>42$.
\end{abstract}

\section{Introduction}

Dissections of regular polygons appear in several popular puzzles (see
\cite{F}). Some of these dissections, such as Langford's dissections of the
regular pentagon \cite{La}, Freese's dissection of the regular octagon
\cite[Figure 17.1]{F}, or K\"ursch\'ak's dissection of the regular $12$-gon
\cite[Figure 2.6.4]{G} consist of triangles of two different shapes.

In this paper we consider dissections of the regular polygons using triangles
of one single shape but not necessarily of the same size. What 
we are interested in is the existence of tilings, independently of the
rearrangement of the pieces (which is the usual motivation for the
puzzles mentioned). We confine our attention to triangles having angles that are
rational multiples of $\pi$. Our aim is to show that if $N$ is large enough,
then there are at most three nonsimilar triangles $T$ in this class
such that the regular $N$-gon can be dissected into similar copies of  $T$.

\subsection{Main results}
By a dissection (or tiling) of a polygon $A$ we mean a decomposition of $A$
into finitely many nonoverlapping polygons. No other conditions are imposed on
the tilings. In particular, it is allowed that two pieces have a common boundary
point, but do not have a common side. We say that a triangle $T$ tiles a
polygon $A$, if $A$ can be dissected into finitely many nonoverlapping triangles
similar to $T$. Our main result is the following.

\begin{theorem}\label{t1}
Suppose that a triangle with angles $\al ,\be ,\ga$ tiles the regular $N$-gon,
where $N\ge 25$ and $N\ne 30, 42$. If $\al ,\be ,\ga$
are rational multiples of $\pi$, then, after a suitable permutation of $\al ,
\be , \ga$, one of the following statements is true:
\begin{enumerate}[{\rm (i)}]
\item $\al =\be =(N-2)\pi /(2N)$ and $\ga =2\pi /N$,
\item $\al =(N-2)\pi /(2N)$, $\be =\pi /N$ and $\ga =\pi /2$, or
\item $\al =(N-2)\pi /N$ and $\be =\ga =\pi /N$.
\end{enumerate}
\end{theorem}

Let $R_N$ and $\de _N$ denote the regular $N$-gon and its angle; that is,
let $\de _N =(N-2)\pi /N$. Connecting the center of $R_N$ with the vertices of
$R_N$ we obtain a dissection of $R_N$ into $N$ congruent isosceles triangles
with angles listed in (i). Bisecting each of these triangles into two right
angled triangles, we get a dissection of $R_N$ into $2N$ congruent triangles
with angles listed in (ii).

\begin{figure}
\centering
\includegraphics[width=2.5in]{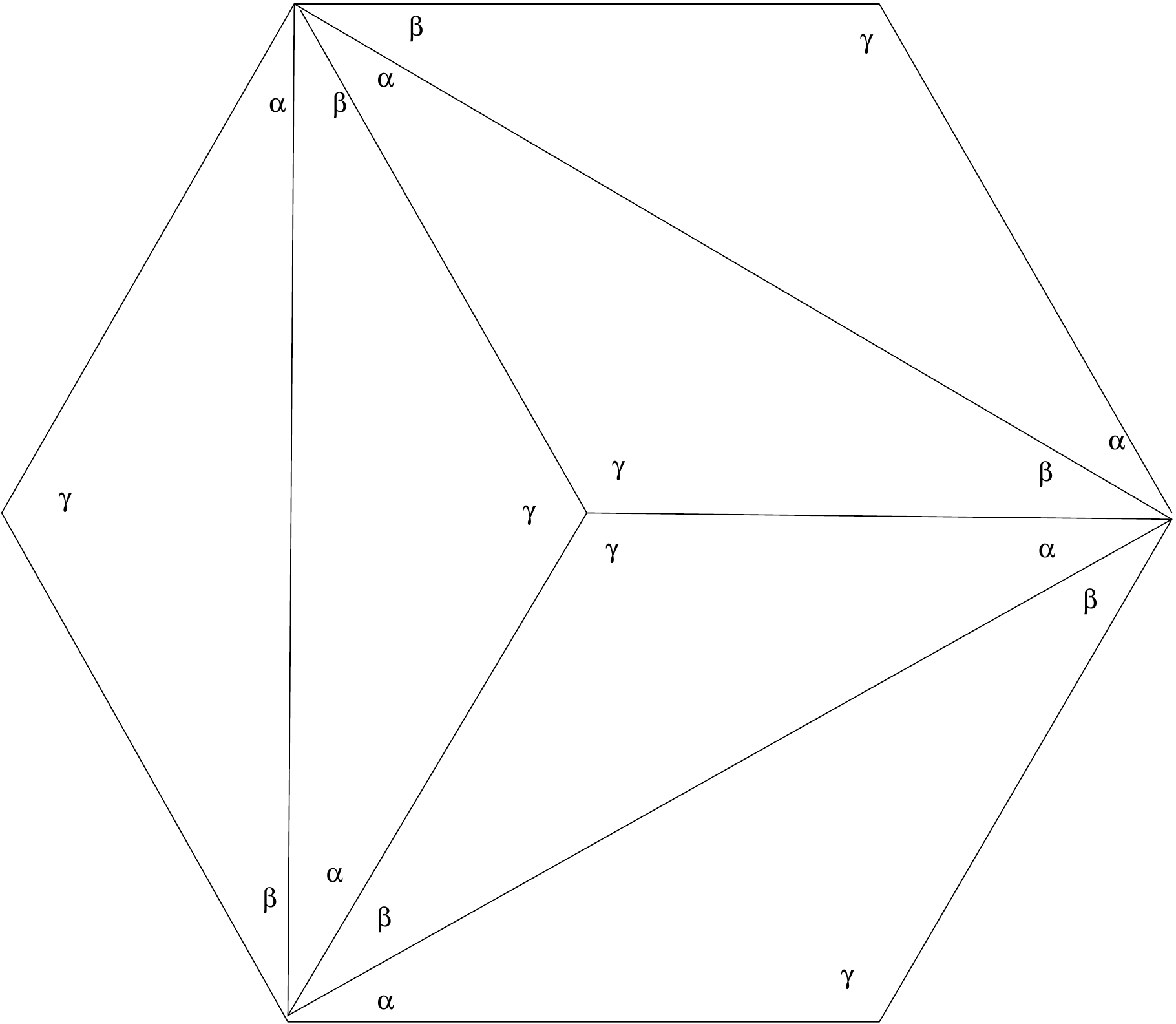}
\caption{a (regular) tiling of $R_6$}
\label{fig1}
\end{figure}

Thus the triangles with angles listed in (i) and (ii) tile $R_N$,
even with congruent copies. This is also true for the triangle
with angles listed in (iii) if $N=3,4$ or $6$. (As for $N=6$, see Figure
\ref{fig1}.) If $N$ is different from $3,4$ or $6$, then dissections of $R_N$
with {\it congruent}
copies of a triangle with angles $\al =\de _N$ and $\be =\ga =\pi /N$ do not
exist (see \cite[Lemma 3.5]{L3}). It is not clear, however, if $R_N$ can be
dissected into {\it similar} triangles of angles $\al =\de _N$ and $\be =\ga =
\pi /N$ for every $N$. In a forthcoming paper \cite{L4} we prove that such
tilings exist for $N=5$ and $N=8$.

\begin{figure}
\centering
\includegraphics[width=2.5in]{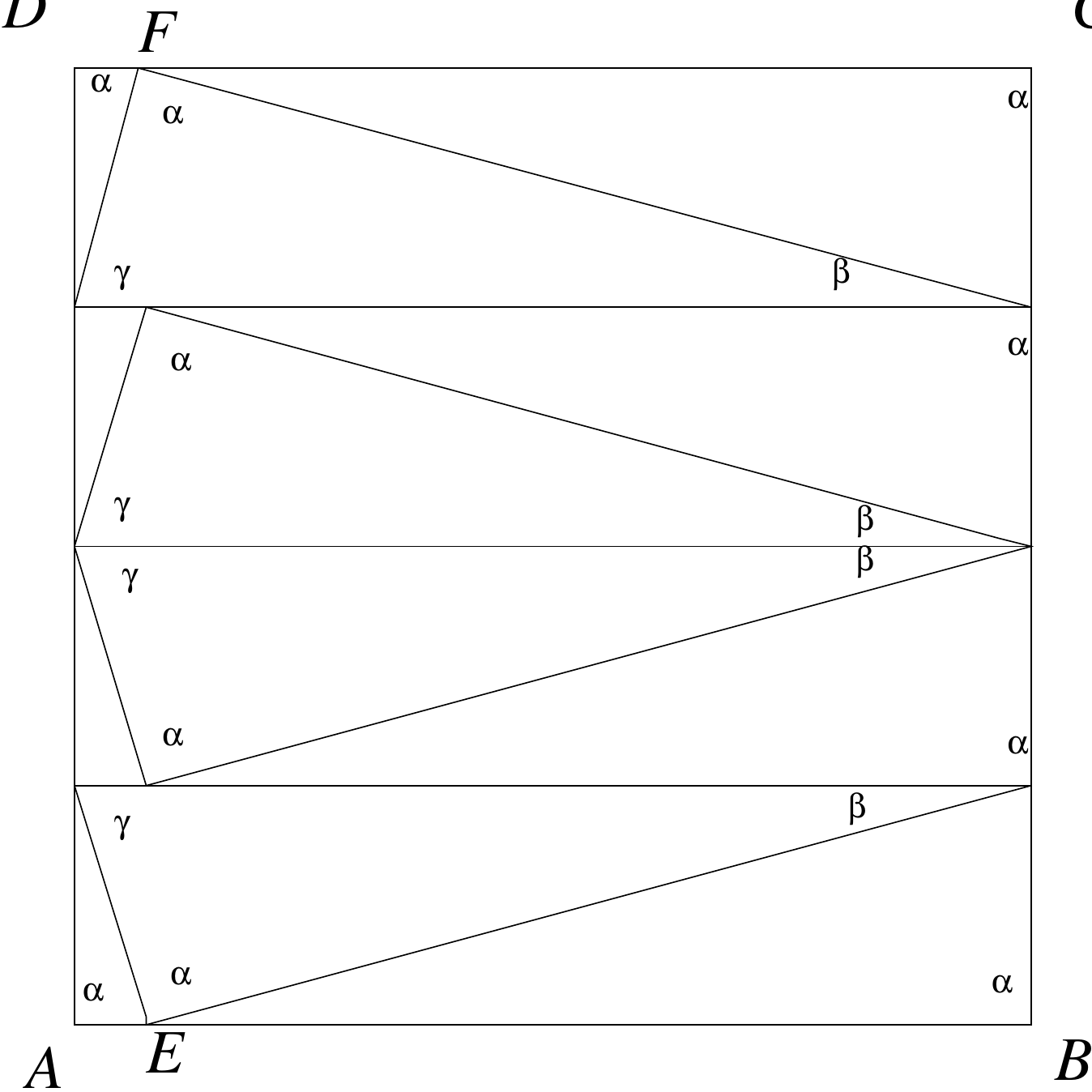}
\caption{a tiling of the square with $\al =\de _4$}
\label{fig2}
\end{figure}

Theorem \ref{t1} will be proved through the following results.
In each of these theorems we assume that {\it a tiling of $R_N$ with triangles
of angles $\al ,\be ,\ga$ is given, where $\al ,\be ,\ga$ are rational multiples
of $\pi$. If the number of angles $\al ,\be , \ga$ meeting at the vertex $V_j$
of $R_N$ is $p_j ,q_j ,r_j$, then we call $p_j \al +q_j \be +r_j \ga =\de _N$
the equation at the vertex $V_j$} $(1\le j\le N)$.
\begin{theorem}\label{l0}
If $N\ne 6$, then we have $p_j +q_j +r_j \le 2$ for every
$j=1\stb N$; that is, each angle of $R_N$ is packed with at most two tiles.
\end{theorem}
Note that the statement of Theorem \ref{l0} is not true for $N=6$, as Figure
\ref{fig1} shows.
\begin{theorem}\label{l1}
Suppose $N>6$. Then the equations at the vertices $V_1 \stb V_N$ are
the same. More precisely, after a suitable permutation of $\al ,\be ,\ga$,
one of the following is true:
\begin{enumerate}[{\rm (i)}]
\item The equation at every vertex $V_j$ is $\al =\de _N$.
\item The equation at every vertex $V_j$ is $\al +\be =\de _N$.
\item The equation at every vertex $V_j$ is $2\al =\de _N$.
\end{enumerate}
\end{theorem}
As Figure \ref{fig1} shows, the statement of Theorem \ref{l1} is not true for
$N=6$. 
\begin{theorem}\label{l2}
Suppose $N>5$. If the equation at every vertex $V_j$ is $\al =\de _N$,
then we have $\be =\ga =\pi /N$.
\end{theorem}
The statement of Theorem \ref{l2} is not true for $N=4$. Figure \ref{fig2}
shows a tiling of the square $ABCD$ with $12$ right triangles
of angles $\al =\pi /2$, $\be =\pi /12$ and $\ga =5\pi /12$.
If the side length of the square is $4$ then we have $\ol{AE}=\ol{DF}=2-\sqrt 3$
and $\ol{EB}=\ol{FC}=2+\sqrt 3$. Note that in this tiling $\al =\de _4$ at each
vertex of the square but $\be \ne \ga$.
We do not know if the statement of Theorem \ref{l2} is true for $N=5$.
\begin{theorem}\label{l3}
Suppose $N>10$. If the equations at the vertices $V_j$ are $\al +\be  =\de _N$,
then we have $\al =\be =\de _N /2$ and $\ga =2\pi /N$.
\end{theorem}
\begin{theorem}\label{l4}
Suppose $N\ge 25$ and $N\ne 30, 42$. If the equations at the vertices $V_j$
are $2\al  =\de _N$, then we have either $\al =\ga =\de _N /2$ and $\be =
2\pi /N$, or $\al =\de _N /2$, $\be =\pi /N$ and $\ga =\pi /2$.
\end{theorem}
It is clear that Theorem \ref{t1} follows from Theorems \ref{l1}-\ref{l4}.
As for the sharpness of the bounds appearing in Theorems \ref{l3} and \ref{l4}
we refer to Remark \ref{r} below.

\subsection{Regular and irregular tilings}
A tiling into similar triangles is called {\it regular}, if the pieces have
two angles, $\al$ and $\be$, such that at each vertex $V$ of any of the tiles,
the number of tiles having angle $\al$ at $V$ is the same as
the number of tiles having angle $\be$ at $V$.
Otherwise the tiling is {\it irregular}. It is known that the number of
triangles that tile a given polygon irregularly is always finite
(see \cite[Theorem 4]{L2}). On the other hand, for every $N\ge 3$ there
are infinitely many triangles that tile the regular $N$-gon regularly (see
\cite[Theorem 2]{L2}).

The problem of listing all triangles that tile a given polygon is difficult;
it is unsolved even for the regular triangle. In fact, the problem is solved
only for the square (see \cite{L} and \cite{Sz}). (See also \cite{L3}, where
the tilings of convex
polygons with congruent triangles are considered.) As for irregular tilings of
$R_N$ $(N>10)$, we have the following corollary of Theorems \ref{l1}-\ref{l3}.
\begin{theorem}\label{t2}
Suppose a triangle $T$ with angles $\al ,\be , \ga$ tiles $R_N$, where $N>10$.
Then there is an irregular tiling of $R_N$ with pieces similar to $T$ if and
only if $\al ,\be , \ga$ are rational multiples of $\pi$.
\end{theorem}
\proof
Suppose there is an irregular tiling of $R_N$ with pieces similar to $T$.
Let $V_1 \stb V_M$ denote the vertices of the tiles, where $M\ge N$ and $V_1
\stb V_N$ are the vertices of $R_N$. If the number of angles $\al ,\be , \ga$
meeting at $V_j$ is $p_j ,q_j ,r_j$, respectively, then we have
$p_j \al +q_j \be +r_j \ga =\si _j $, where $\si _j =\de _N$ if $j=1\stb N$,
and $\si _j$ equals $\pi$ or $2\pi$ if $N<j\le M$. If the tiling is irregular,
then, by \cite[Lemma 10]{L2}, there are indices $i<j$ such that the determinant
$D_{ij}=\left| 
\begin{matrix}
1&1&1\cr  p_i & q_i & r_i \cr
p_j & q_j &r_j 
\end{matrix}
 \right|$ is nonzero. Then the corresponding system of equations
\begin{align*}
& \phantom{p_i} \al +\phantom{q_i} \be +\phantom{r_i} \ga =\pi \\
& p_i \al +q_i \be +r_i \ga =\si _i \\
& p_j \al +q_j \be +r_j \ga =\si _j 
\end{align*}
determines $\al ,\be ,\ga$. Applying Cramer's rule, we find that $\al ,\be ,\ga$
are rational multiples of $\pi$. 

Next let $\al ,\be , \ga$ be rational multiples of $\pi$. Since $N>10$, one of
(i), (ii) and (iii) of Theorem \ref{l1} holds. If (i) or (ii) holds, then it
follows from Theorems \ref{l2} and \ref{l3} that $T$ is isosceles. Suppose
$\al =\be$, and consider a tiling of $R_N$ with pieces similar to $T$.
If the tiling is irregular, we are done. If, however, it is regular, then
changing the labels $\al$ and $\be$ in one of the pieces we obtain an irregular
tiling.

Now suppose that (iii) of Theorem \ref{l1} holds. We prove that in this case
every tiling with similar copies of $T$ must be irregular.
Suppose this is not true, and consider a regular tiling. Since the equation
at each vertex of $R_N$ is $2\al =\de _N$, it follows that
$q_j =r_j $ for every $j$. Then there must be an equation with $p_j <q_j =r_j$,
since in the equations at the vertices we have $p_j >q_j =0$. For this equation
we have
$$(q_j -p_j )(\be +\ga )=(p_j \al +q_j \be +r_j \ga )-p_j (\al +\be +\ga )=
v_j \pi -p_j \pi =(v_j -p_j )\pi ,$$
hence $(q_j -p_j )((1/2)+(1/N))=
v_j -p_j$ and $(q_j -p_j )\cd (N+2)=2(v_j -p_j )N$. Since $q_j -p_j $ is a
positive integer, we have $(N+2) \mid 2(v_j -p_j )N$ and
$N+2 \mid 4(v_j -p_j )$. Now $v_j -p_j$ is positive, since $(q_j -p_j )\cd (N+2)>
0$. Then $0<v_j -p_j \le 2$, $0<4(v_j -p_j )\le 8$, and thus $(N+2) \mid
4(v_j -p_j )$ implies $N\le 6$, which is impossible. \hfill $\square$

\bsk
Comparing Theorem \ref{t2} with Theorem \ref{t1} we obtain the following.
\begin{cor}\label{c}
If $N>42$, then there are at most three triangles
that tile the regular $N$-gon irregularly.
\end{cor}

\subsection{Condition (K) and Condition (E)}
The main tool in the proof of Theorems \ref{l0}-\ref{l4} is the next result.
\begin{lemma}\label{l5}
Suppose $R_N$ can be dissected into finitely many triangles with angles
$\al =(a/n)\pi ,\ \be =(b/n)\pi ,\ \ga  =(c/n)\pi$,
where $a,b,c,n$ are positive integers with $a+b+c=n$.
Let the equation at the vertices of $R_N$ be $p_j \al +q_j \be +r_j \ga =\de _N$
$(j=1\stb N)$.

If $k$ is prime to $n\cd N$ and $\{ k/N\} <1/2$, then we have
\begin{equation}\label{e13}
\left\{ \frac{ka}{n} \right\} + \left\{ \frac{kb}{n} \right\} +
\left\{ \frac{kc}{n} \right\} =1
\end{equation}
and
\begin{equation}\label{e12}
p_j \left\{ \frac{ka}{n} \right\} + q_j \left\{ \frac{kb}{n} \right\} +
r_j \left\{ \frac{kc}{n} \right\} =1-2\left\{ \frac{k}{N} \right\} 
\end{equation}
for every $j=1\stb N$.
\end{lemma}

We say that {\it the angles $\al =(a/n)\pi ,\ \be =(b/n)\pi ,\ \ga =
(c/n)\pi$ satisfy Condition} (K), if the conclusion of the lemma above
holds; that is, if \eqref{e13} and \eqref{e12} hold true for every $k$
such that $\gcd (k, nN)=1$ and $\{ k/N\} <1/2$.
As we shall see in the next section, Condition (K) is deduced from the
properties of conjugate tilings.

If a tiling exists with triangles of angles $\al ,\be , \ga$, then the angles
have to satisfy another necessary condition: there must exist nonnegative
integers $p_j, q_j, r_j$ $(j=1\stb M; \ M\ge N)$ such that
\begin{enumerate}[{\rm (i)}]
\item $p_j \al +q_j \be +r_j \ga =\de _N$ for every $j=1\stb N$, 
\item $p_j \al +q_j \be +r_j \ga$ equals $\pi$ or $2\pi$ for every
$j=N+1\stb M$, and 
\item $\sum_{j=1}^M p_j =\sum_{j=1}^M q_j =\sum_{j=1}^M r_j$.
\end{enumerate}
We say that {\it the angles $\al , \be ,\ga$
satisfy Condition} (E), if there are nonnegative integers $p_j, q_j , r_j$
with these properties.

In the proof of Theorems \ref{l0}-\ref{l4} we only use Condition (K) and
Condition (E) on the angles $\al ,\be ,\ga$. In fact, I am not aware of any
other necessary condition that must be satisfied by the angles of a tiling,
if they are rational multiples of $\pi$. Perhaps it would be hasty to conjecture
that whenever the angles of a triangle satisfy Condition (K) and Condition (E),
then a tiling must exist. Still, it should be remarked that tilings of $R_N$
with triangles of angles $\al =\de _N$ and $\be =\ga =\pi /N$ were found at
least for the regular pentagon and octagon \cite{L4}. In this context
I also mention B. Szegedy's remarkable tilings of the square with right
triangles, found ten years after the necessary conditions were established
\cite{Sz}.
\begin{remark}\label{r}
{\rm We do not know if the lower bounds in Theorems \ref{l2}-\ref{l4} are sharp
or not. We show, however, that if we only use Condition (K) and Condition (E),
then these bounds cannot be improved. As for Theorem \ref{l2}, consider the 
triangle $T_1$ with angles
$$(\al ,\be ,\ga )=
\left( \frac{6\pi}{10} , \frac{\pi}{10}, \frac{3\pi}{10} \right) .$$
Then the existence of a tiling of $R_5$ with similar copies of $T_1$ cannot be
disproved by only using Condition (K) and Condition (E). Indeed, suppose
that the equation at each vertex of $R_5$ is $\al =\de _5$. Then Condition (K)
is satisfied. Indeed, the only $k$ with $1<k<10$, $\gcd (k,10)=1$ and $\{ k/5\}
<1/2$ is $k=7$, and it is easy to check that both \eqref{e13} and \eqref{e12}
are satisfied if $(a/n,b/n,c/n)=(6/10,1/10,3/10)$ and $k=7$.
Condition (E) is also satisfied. Indeed, consider the following system of
equations: take $5$ equations $\al =\de _5$, an equation $\be +3\ga =\pi$ and
an equation $4\be +2\ga =\pi$.

As for Theorem \ref{l3}, consider the triangle $T_2$ with angles
$$(\al ,\be ,\ga ) =
\left( \frac{7\pi}{10} , \frac{\pi}{10}, \frac{2\pi}{10} \right) .$$
Then the existence of a tiling of $R_{10}$ with similar copies of $T_2$ cannot be
disproved by only using Condition (K) and Condition (E). Suppose that the
equation at each vertex of $R_{10}$ is $\al +\be =\de _{10}$. Then Condition (K)
is satisfied. Indeed, the only $k$ with $1<k<10$, $\gcd (k,10)=1$ and
$\{ k/10\} <1/2$ is $k=3$, and it is easy to check that both \eqref{e13} and
\eqref{e12} are satisfied if $(a/n,b/n,c/n)=(7/10,1/10,2/10)$ and $k=3$.
Condition (E) is also satisfied: take $10$ equations $\al +\be =\de _{10}$
and an equation $10\ga =2\pi$. 

In the case of Theorem \ref{l4}, consider the triangle $T_3$ with angles
$$(\al ,\be ,\ga ) =\left( \frac{20\pi}{42} , \frac{10\pi}{42}, \frac{12\pi}{42}
\right) ,$$
and let the equation at each vertex of $R_{42}$ be
$2\al =\de _{42}$. Then Condition (K) is satisfied. Indeed, if $1<k<42$,
$\gcd (k,42)=1$ and $\{ k/42\} <1/2$, then $k$ is one of $5,11,13,17,19$.
It is easy to check that both \eqref{e13} and \eqref{e12}
are satisfied if $(a/n,b/n,c/n)=(20/42,10/42,12/42)$ and if $k$ is any of these
values. Condition (E) is also satisfied: take $42$ equations $2\al =\de _{42}$,
$8$ equations $7\ga =2\pi$ and $28$ equations $3\be +\ga =\pi$.

Similarly, if $N=30$, then the triple
$$\left( \frac{14\pi}{30}, \frac{6\pi}{30}, \frac{10\pi}{30} \right)$$
satisfies both Condition (K) and Condition (E). As for the latter,
take $30$ equations $2\al =\de _{30}$,
$20$ equations $3\ga =\pi$ and $12$ equations $5\be  =\pi$. 
}
\end{remark}

\subsection{Further lemmas}
Since Condition (K) is of arithmetical nature, it can be expected that
in the arguments involving Condition (K) we need some facts of elementary
number theory. These facts are collected in the next lemmas. Their proofs, being
independent of the rest of the paper, are postponed to the last three sections.
\begin{lemma}\label{lL7}
Let $a,n,N,N'$ be positive integers such that $\gcd(a,n)=1$ and $\gcd (N,N')=1$.
Then one of the following statements is true.
\begin{enumerate}[{\rm (i)}]
\item There exists an integer $k$ such that $\gcd (k,nN)=1$, $k\equiv N'$
  (mod $N$), and $\{ ka/n\} \ge 1/3$.
\item $N$ is odd and $n\mid 2N$.
\item $N$ is even and $n\mid N$.
\end{enumerate}
\end{lemma}

\begin{lemma}\label{lL11}
Let $a,b,n,N$ be positive integers and $p,q$ be nonnegative integers such that
$a+b<n$, $N\ge 3$, $N\ne 6$, and 
\begin{equation}\label{eL11}
p  \left\{ \frac{ka}{n} \right\} +q  \left\{ \frac{kb}{n} \right\} = 
1-2\left\{ \frac{k}{N}  \right\}
\end{equation}
for every integer $k$ satisfying $\gcd (k,nN)=1$ and $\{ k/N\} <1/2$.
Then we have $p+q\le 2$.
\end{lemma}
Note that Theorem \ref{l0} is an immediate consequence of Lemmas \ref{l5} and
\ref{lL11}.
\begin{lemma}\label{lL1}
\begin{enumerate}[{\rm (i)}]
\item For every even integer $N\ge 26$ there are integers $k, k'$ such that
$N/4 < k,k'<N/2$, $\gcd (k,N)=\gcd (k',N)=1$, $k\equiv 1$ (mod $4$), and
$k'\equiv 3$ (mod $4$).
\item For every $N\ge 43$ there exists an integer $k$ such that
$N/6 < k<N/4$ and $\gcd (k,2N)=1$.
\end{enumerate}
\end{lemma}

The following simple observation will be used frequently.
\begin{proposition}
Let $u,v,n$ be nonzero integers. If $\gcd(u,v)=1$, then there exists an integer
$j$ such that $u+jv$ is prime to $n$.  
\end{proposition}\label{p}
\proof
Let $j$ be the product of those primes that divide $n$ but does not divide
$u$. (We put $j=1$ if there is no such prime.) Then every prime divisor of $n$
divides exactly one of $u$ and $jv$, and thus $\gcd (u+jv,n)=1$. \hfill $\square$

\bsk
The paper is organized as follows. In the next five sections we prove Lemma
\ref{l5} and Theorems \ref{l1}-\ref{l4}, in this order. Then we prove Lemmas
\ref{lL7}-\ref{lL1} in Sections 7-9.

\section{Proof of Lemma \ref{l5}}
Let the vertices of $R_N$ be the $N^{\rm th}$ roots of unity; that is, let
$V_j =e^{2\pi ji/N}$ for every $j=0\stb N-1$. First we assume that $4N\mid n$.
Let $\zeta$ denote the first $n^{\rm th}$ root of unity, and let $F$ denote the
field of real elements of the cyclotomic field $\qq (\zeta )$. Then the
coordinates of the vertices of $R_N$ belong to $F$, 
since $\cos 2j \pi /N =(\zeta ^{nj/N} +\zeta ^{-nj/N})/2$ and 
$\sin 2j \pi /N =(\zeta ^{nj/N} -\zeta ^{-nj/N})/(2\zeta ^{n/4})$ for every
integer $j$. Also, $\cot \al ,\cot \be ,\cot \ga$ belong to $F,$ since
$$\cot \frac{j}{n}\pi =\frac{e^{(j/n)\pi i} +e^{-(j/n)\pi i}}{e^{(j/n)\pi i} -
e^{-(j/n)\pi i}} \cd \zeta ^{n/4}=\frac{\zeta ^j +1}{\zeta ^j -1} \cd \zeta ^{n/4}$$
for every $j$. Let $\De _1 \stb \De _t$ be the tiles of the dissection.
By Theorem 1 of \cite{L}, the coordinates of the vertices of the triangles
$\De _j$ belong to $F$. 

Let $k$ be an integer prime to $n$, and let $\ph \colon \qq (\zeta) \to \cc$
be the isomorphism of $\qq (\zeta )$ satisfying $\ph (\zeta )=\zeta ^k$. Then
$\ph$ commutes with complex conjugation, and thus 
$\phi $ restricted to $F$ is also an isomorphism. It is easy to check that
$$\ph \left( \cot \frac{j}{n}\pi \right) =(-1)^{(k-1)/2 } \cot \frac{kj}{n}\pi $$
for every integer $j$. We define $\Ph (x,y)=(\ph (x),\ph (y))$ for every $x,y\in
F$. Then $\Ph$ is a collineation defined on $F\times F$. In particular, $\Ph$ is
defined on the set of vertices of the tiles $\De _j$ $(j=1\stb t)$. We denote by
$\De '_j$ the triangle with vertices  $\Ph \left( V_{j,1} \right) ,
\Ph \left( V_{j,2} \right) ,\Ph \left( V_{j,3} \right) ,$
where $V_{j,1},V_{j,2},V_{j,3},$ are the vertices of $\De _j $.

Let $\ep _j =1$ if $\Ph$ does not change the orientation of $\De _j ,$ and let 
$\ep _j =-1$ otherwise. If the angles of $\De ' _j$ are $\al `_j ,\be `_j ,
\ga `_j ,$ then, by Lemma 6 of \cite{L}, we have
$$\cot \al '_j =\ep _j \cd \ph (\cot \al )=\ep _j \cd (-1)^{(k-1)/2} \cd \cot \frac{ka}{n}\pi$$
and, similarly,
$$\cot \be '_j =\ep _j \cd (-1)^{(k-1)/2} \cd \cot \frac{kb}{n}\pi ,\qquad
\cot \ga '_j =\ep _j \cd (-1)^{(k-1)/2} \cd \cot \frac{kc}{n}\pi .$$
Note that at least two of the numbers $\cot \al '_j ,\cot \be '_j ,\cot \ga '_j$
are positive for every $j$. Since the integers $a,b,c,n,k$ are fixed, this
implies that the value of $\ep _j$ is the same for every $j=1\stb t$. Therefore,
the orientation of the triangles $\De '_j$ is the same, and the angles of each
$\De '_j$ are
\begin{equation}\label{e10}
\al '=\left\{ \frac{ka}{n} \right\} \pi , \quad  \be '=\left\{ \frac{kb}{n} \right\} \pi , \quad 
\ga '=\left\{ \frac{kc}{n} \right\} \pi
\end{equation}
if $\ep \cd (-1)^{(k-1)/2}=1,$ and
\begin{equation}\label{e11}
\al '=\left( 1-\left\{ \frac{ka}{n} \right\} \right) \pi , \quad  
\be '=\left( 1- \left\{ \frac{kb}{n} \right\} \right)  \quad
\ga '=\left( 1- \left\{ \frac{kc}{n} \right\} \right) \pi
\end{equation}
if $\ep \cd (-1)^{(k-1)/2}=-1,$ where $\ep$ is the common value of $\ep _j \ (j=1\stb t)$.

Note that by $4\mid n$ we have $i=\zeta ^{n/4} \in \qq (\zeta )$ 
and $\ph (i)=\zeta ^{kn/4}=(-1)^{(k-1)/2} \cd i$.
If we identify $\sik$ with $\cc$ then we find that for every
$z=x+iy \in \qq (\zeta )$ we have $\Ph (z)= \ph (x)+i\ph (y)=\ph (z)$ if
$(-1)^{(k-1)/2}=1$, and $\Ph (z)= \ol{\ph (z)}$ if $(-1)^{(k-1)/2}=-1$. 

Clearly, $\Ph (V_1 )\stb \Ph (V_N )$ are the vertices of
a star polygon $R'_N$. By the previous observation, the order of the vertices of
$R'_N$ are $1,\zeta ^{kn/N} \stb \zeta ^{(N-1)kn/N}$ or $1,\zeta ^{-kn/N} \stb
\zeta ^{-(N-1)kn/N}$ depending on the sign of $(-1)^{(k-1)/2}$. 

Suppose $\{ k/N \} <1/2$. Then the angles of $R'_N$ at the vertices equals
$(1-2\{ k/N\} )\pi ,$ and the orientation of $R'_N$ is positive or negative
according to the sign of $(-1)^{(k-1)/2}$.

Let $w(x;P)$ denote the winding number of a closed polygon $P$ at a point
$x\notin P;$ that is, let $w(x;P)=(1/(2\pi i)) \int_P dz/(z-x)$. Since the
boundary $\del R'_N$ of $R'_N$ as an oriented cycle equals the sum of the
boundaries $\del \De '_j ,$ we have
$$w(x;\del R'_N )=\sum_{j=1}^t w(x;\del \De '_j ).$$
If $x$ does not belong to the boundaries of $\De ' _j,$ then we have either 
$w(x;\del \De '_j ) =\ep$ or $w(x;\del \De '_j ) =0$ for every $j$. Therefore,
if $w(x;R'_N)=\pm 1,$ then $x$ belongs to exactly one of the
triangles $\De '_j$. Now, for each vertex $V'_j$ $(j=1\stb N)$ there is angular
domain $D_j$ of angle $(1-2\{ k/N\} )\pi$ and there is a neighbourhood $U_j$ of
$V_j$ such that $w(x;R'_N)=(-1)^{(k-1)/2}$ if $x\in U_j \cap D_j$, and 
$w(x;R'_N)=0$ if $x\in U_j \se D_j$. This implies that $\ep =(-1)^{(k-1)/2}$,
the triangles having a vertex at $V'_j$ are nonoverlapping, and their union in
$U_j$ equals $U_j \cap D_j $. Therefore, the angles $\al ',\be ', \ga '$ are
given by \eqref{e10}, and thus \eqref{e13} and \eqref{e12} hold.
This proves the theorem in the case when $4N\mid n$.

In the general case we put
$m=4Nn$. Then we have $\al =(4Na/m)\pi$, $\be =(4Nb/m)\pi$, $\ga =(4Nc/m)\pi$.

Let $k$ be prime to $n\cd N,$ and suppose $\{ k/N\} <1/2$. Then $k+snN$ is
prime to $m$ for a suitable $s$ by Proposition \ref{p}. Since $\{ (k+snN)/N\}
=\{ k/N\} <1/2 ,$ and $\{ (k+snN)\cd 4Na/m\} =\{ ka/n\}$ etc., we obtain
\eqref{e13} and \eqref{e12}. \hfill $\square$

\section{Proof of Theorem \ref{l1}}

In the next two sections we write $\de$ for $\de _N$. By Theorem \ref{l0}, the
equation at each vertex $V_j$ equals one of $\al =\de$, $\be =\de$, $\ga =\de$,
$\al +\be =\de$, $\al +\ga =\de$, $\be +\ga =\de$, $2\al =\de$, $2\be =\de$,
$2\ga =\de$.

First suppose that $\al =\de$ is one of the equations. If $\be =\de$ is another,
then $\al +\be <\pi$ gives $2\de <\pi$, $2(N-2)/N <1$ and $N<4$, which is
impossible. We have the same conclusion if $\ga =\de$.

It is clear that $\al +\be =\de$ or $\al +\ga =\de$ is impossible. If
$\be +\ga =\de$, then $2\de =\al +\be +\ga =\pi$, $\de =\pi /2$ and $N=4$,
which is impossible.

Clearly, $2\al =\de$ is impossible. If $2\be =\de$, then $\al +\be <\pi$ gives
$\al +\be =3\de /2 =(3\pi /2) -(3\pi /N ) <\pi$ and $N<6$, which is impossible.
We have the same conclusion if $2\ga =\de$.
We find that if $\al =\de$ is one of the equations, then each of the equations
is $\al =\de$, and we have (i).

Therefore, we may assume that the equation at each vertex $V_j$ equals one of
$\al +\be =\de$ etc., $2\al =\de$ etc.

Suppose that $\al +\be =\de$ is one of the equations. If $\al +\ga =\de$ is
another, then $\be =\ga$, $\al =\pi -2\be$, $\de =\pi -\be$,
$\be =\ga =2\pi /N$, $\al =(N-4)\pi /N$. Let
$$k=
\begin{cases}
(N-1)/2 &\text{if $N$ is odd},\\
(N/2)-1 &\text{if $4\mid N$},\\
(N/2)-2 &\text{if $N\equiv 2$ (mod $4$)}.
\end{cases}
$$
Then $\gcd (k,N)=1$ and $0<k<N/2$. By Lemma \ref{l5}, this implies that
\eqref{e13} holds, hence $4k/N=\{ 2k/N\} +\{ 2k/N\} <1$ and $k<N/4$.
If $N$ is odd, then this implies $(N-1)/2<N/4$, which is impossible.
If $4\mid N$, then $(N/2)-1<N/4$ is also impossible. If $N\equiv 2$ (mod $4$),
then we get $(N/2)-2<N/4$, $N<8$, $N=6$, which is excluded.
We have the same conclusion if $\be +\ga =\de$.

If $2\ga =\de$ is another equation, then $\pi =\al +\be +\ga =3\de /2$,
$3(N-2)/(2N)=1$ and $N=6$, which is impossible.

We find that if $\al +\be =\de$ is one of the equations, then either 
each of the equations is $\al +\be =\de$; that is, (ii) holds, or each of
the other equations is one of $\al +\be =\de$, $2\al =\de$ and $2\be =\de$,
and at least one of $2\al =\de$ and $2\be =\de$ must occur. Then we have
$\al =\be =\de /2$, and
the tiles are isosceles. It is easy to check that in this case we can exchange
the labels of the angles $\al$ and $\be$ in some of the tiles such that each
equation at the vertices becomes $\al +\be =\de$, and thus (ii) holds.

Therefore, we may assume that the equation at each vertex $V_j$ equals one of
$2\al =\de$, $2\be =\de$ and $2\ga =\de$. If all of these equations occur,
then $\al =\be =\ga =\pi /3$, $\de =2\pi /3$ and $N=6$, which is excluded.

If two of them, say $2\al =\de$ and $2\be =\de$ occur, then we have $\al =
\be =\de /2$, and the tiles are isosceles. Then, as above, we can exchange
the labels of the angles $\al$ and $\be$ in some of the tiles such that each
equation at the vertices becomes $\al +\be =\de$, and thus (ii) holds.

Finally, if only $2\al =\de$ occurs, then we have (iii). \hfill $\square$

\section{Proof of Theorem \ref{l2}}
We have $\al =\de$ and $\be +\ga =2\pi /N$. If $\be =\ga$, then we have
$\be =\ga =\pi /N$, and we are done. Therefore, we may assume $\ga >
\be$ by symmetry.

Let $\al =(a/n)\pi$, $\be =(b/n)\pi$, $\ga =(c/n)\pi$, where $a,b,c,n$ are
positive integers such that $a+b+c=n$.
Let $b/n =b_2 /n_2$ and $c/n = c_3 /n_3$, where $\gcd (b_2 ,n_2 )=
\gcd (c_3 ,n_3 )=1$.

First we suppose $n_3 \mid N$. Then
\begin{equation}\label{e2}
\frac{b_2}{n_2} =1- \frac{\ga}{\pi} -\frac{\al}{\pi} =
1-\frac{c_3}{n_3}-\frac{N-2}{N}
\end{equation}
gives $n_2 \mid N$. Thus $b_2 /n_2 \ge 1/N$ and
$c_3 /n_3 \ge 1/N$. Since $(b_2 /n_2 )+(c_3 /n_3 )=2/N$ and $c_3 /n_3 > b_2 /n_2$,
this is impossible.

Next suppose that $N$ is odd and $n_3 \mid 2N$. Then \eqref{e2} gives $n_2 \mid
2N$. Now we have $(b_2 /n_2 )+(c_3 /n_3 )=2/N=4/(2N)$, and thus we have
$c_3 /n_3 =3/(2N)$ and $b_2 /n_2 =1/(2N)$ by $c_3 /n_3 > b_2 /n_2$.

Since $k=N+2$ is prime to $2N$ and $\{ k/N\} =2/N<1/2$, it follows from
\eqref{e13} that
$$\left\{ \frac{3(N+2)}{2N} \right\} + \left\{ \frac{N+2}{2N} \right\} <1,$$
which is absurd.

Therefore, we may assume $n_3 \ {\mid \!\!\! \not} \quad N$ and, if $N$ is odd,
then $n_3 \ {\mid \!\!\! \not} \quad 2N$. By Lemma \ref{lL7}, this
implies that there is a $k$ prime to $n_3 N$ and such that $k\equiv 1$ (mod $N$)
and $\{ kc_3  /n_3 \} \ge 1/3$. Replacing $k$ by $k+jn_3 N$ with a suitable $j$,
we may assume that $k$ is prime to $nN$.

Then \eqref{e13} gives $\{ ka/n\} +\{ kc/n\} <1$. Since 
$\{ ka/n\} =\{ k(N-2)/N\}  =\{ (N-2)/N\} =(N-2)/N$ and $\{ kc/n\} \ge 1/3$,
we have $(N-2)/N <2/3$ and $N<6$, which is impossible. \hfill $\square$

\section{Proof of Theorem \ref{l3}}
We put
$$N'=
\begin{cases}
(N-1)/2 &\text{if $N$ is odd},\\
(N/2)-1 &\text{if $4\mid N$},\\
(N/2)-2 &\text{if $N\equiv 2$ (mod $4$)}.
\end{cases}
$$
Then $\gcd (N,N')=1$ and $\{ N'/N\} <1/2$.

Let $\al = a_1 \pi /n_1$, where $\gcd (a_1 ,n_1 )=1$. By Lemma \ref{lL7},
at least one of the following statements is true: (i) there exists a $k$ such
that $k\equiv N'$ (mod $N$), $\gcd (k,n_1 N)=1$ and $\{ ka_1 /n_1 \} \ge 1/3$,
(ii) $N$ is odd and $n_1 \mid 2N$, and (iii) $n_1 \mid N$.

If (i) holds then we may assume that $k$ also satisfies $\gcd (k,nN)=1$.
Indeed, if $k$ satisfies the conditions of (i), then so does $k+jn_1 N$
for every $j$. Replacing $k$ by $k+jn_1 N$ with a suitable $j$, we find that 
$\gcd (k,nN)=1$ will also hold.
Therefore, by \eqref{e12} of Lemma \ref{l5}, we obtain
$$\frac{1}{3} \le \left\{ \frac{ka_1}{n_1}\right\} =
\left\{ \frac{ka}{n}\right\} <1-2\cd \left\{ \frac{k}{N}\right\}  =
1-2\cd \frac{N'}{N}$$
and $N'<N/3$, which is impossible by $N> 10$.

Next suppose that $N$ is odd and $n_1 \mid 2N$. Let $b/n=b_2 /n_2$, where
$\gcd (b_2 ,n_2 )=1$. Since $\be =(\al +\be )-\al = (N-2)\pi /N  -\al$,
we have 
\begin{equation}\label{er16}
\frac{b_2}{n_2} =\frac{N-2}{N} -\frac{a_1}{n_1} ,
\end{equation}
and thus $n_2 \mid 2N$. Then $c_3 /n_3 =1-(a_1 /n_2 )-(b_2 /n_2 )$ gives
$n_3 \mid 2N$. Therefore, we may assume $n=2N$.

We put $k=N'$ if $N'$ is odd, and $k=N'+N$ if $N'$ is even. Then $\gcd (k,2N)=1$
and $\{ k/N\} <1/2$, and thus \eqref{e12} of Lemma \ref{l5} gives
\begin{equation}\label{er15}
\begin{split}
\left\{ \frac{ka}{2N}\right\} + \left\{ \frac{kb}{2N}\right\} &=
1- 2\cd \left\{ \frac{k}{N}\right\} =1-2\cd \left\{ \frac{N'}{N}\right\}  =\\
&= 1-2\cd \frac{N'}{N} =1-2\cd \frac{(N-1)/2}{N} =\frac{1}{N}.
\end{split}
\end{equation}
Since $\{ ka/(2N)\}$ and $\{ kb/(2N)\}$ are positive integer multiples of
$1/(2N)$, \eqref{er15} gives $\{ ka/(2N)\} = \{ kb/(2N)\} =1/(2N)$. Then
$ka\equiv kb \equiv 1$ (mod $2N$). By $\gcd (k,2N)=1$ this implies
$a\equiv b$ (mod $2N$), $a=b$, and $\al =\be =((1/2)-(1/N))\pi$.
That is, the statement of the theorem is true in this case.

Finally, suppose that $2\mid N$ and $n_1 \mid N$. Then we may assume $n=N$
by \eqref{er16}. Then \eqref{e12} of Lemma \ref{l5} gives
\begin{equation}\label{er17}
\left\{ \frac{N'a}{N}\right\} + \left\{ \frac{N'b}{N}\right\} =
1- 2\cd \left\{ \frac{N'}{N}\right\} = 1-2\cd \frac{N'}{N}.
\end{equation}
The value of $N'/N$ is $(1/2)-(1/N)$ if $4\mid N$, and $(1/2)-(2/N)$ if
$N\equiv 2$ (mod $4$). Thus $1-2\cd (N'/N)$ equals either $2/N$ or $4/N$.
Since $\{ N'a/N\}$ and $\{ N'b/N\}$ are positive integer multiples of
$1/N$, we have the following possibilities:
$\{ N'a/N\} = \{ N'b/N\} =1/N$, $\{ N'a/N\} = \{ N'b/N\} =2/N$, or $\{ \{ N'a/N\} , \{ N'b/N\} \} = \{ 1/N, 3/N\}$. In the third case we may assume, by
symmetry, that $\{ N'a/N\} =1/N$.

In the first two cases we have $N'a\equiv N'b$ (mod $N$), 
$a\equiv b$ (mod $N$), $a=b$, $\al =\be =((1/2)-(1/N))\pi$, and we are done.

Therefore, we may assume that $N\equiv 2$ (mod $4$) and 
$\{ N'a/N\} =1/N$; that is, $N'a\equiv 1$ (mod $N$). Since $N$ is even, $a$
must be odd. Now $N/2$ is odd either, and thus $(N/2)a\equiv N/2$ (mod $N$).
Then, by $N'=(N/2)-2$ we obtain
$$1\equiv N'a =\left(\frac{N}{2} -2\right) a \equiv \frac{N}{2} -2a \quad
\text{(mod $N$)},$$
and $2a+1\equiv N/2$ (mod $N$). Since $0<a<N$, we have either $2a+1=N/2$;
that is, $a=(N/4)-(1/2)$, or $2a+1=3N/2$; that is, $a=(3N/4)-(1/2)$.

Suppose $a=(N/4)-(1/2)$. By (i) of Lemma \ref{lL1}, if $N\ge 26$,
then there is a $k$ such that $N/4<k<N/2$, $\gcd (k,N)=1$ and $k\equiv 3$
(mod $4$). Then
$$\left\{ \frac{ka}{N}\right\} = \left\{ \frac{k}{4} - \frac{k}{2N}\right\} =
\left\{ \frac{3}{4} -\frac{k}{2N} \right\} =\frac{3}{4} -\frac{k}{2N} >\frac{1}{2}.$$
On the other hand, \eqref{e12} of Lemma \ref{l5} gives
$$\left\{ \frac{ka}{N}\right\} <
1- 2\cd \left\{ \frac{k}{N}\right\} <\frac{1}{2},$$
a contradiction.

If $N=14$, then $a=(N/4)-(1/2)=3$ and $b=N-2-a=9$.
In this case $k=3$ is prime to $14$, $3/14<1/2$, but
$$\left\{ \frac{ka}{N}\right\} +\left\{ \frac{kb}{N}\right\} =
\left\{ \frac{9}{14}\right\} +\left\{ \frac{27}{14}\right\} =
\frac{22}{14} >1,$$
a contradiction.

If $N=18$, then $a=4$ and $b=12$. Then $k=7$ is prime to $18$, $7/18<1/2$, but
$$\left\{ \frac{ka}{N}\right\} +\left\{ \frac{kb}{N}\right\} =
\left\{ \frac{28}{18}\right\} +\left\{ \frac{84}{18}\right\}=\frac{22}{18} >1,$$
a contradiction.

If $N=22$, then $a=5$ and $b=15$. Then $k=7$ is prime to $22$, $7/22<1/2$, but
$$\left\{ \frac{ka}{N}\right\} +\left\{ \frac{kb}{N}\right\} =
\left\{ \frac{35}{22}\right\} +\left\{ \frac{105}{22}\right\}=\frac{30}{22} >
1,$$
a contradiction. Therefore, the case $a=(N/4)-(1/2)$ is impossible if $N>10$.

Next, let $a=(3N/4)-(1/2)$. By (i) of Lemma \ref{lL1}, if $N\ge 26$,
then there is a $k$ such that $N/4<k<N/2$, $\gcd (k,N)=1$ and $k\equiv 1$
(mod $4$). Then
$$\left\{ \frac{ka}{N}\right\} = \left\{ \frac{3k}{4} - \frac{k}{2N}\right\} =
\left\{ \frac{3}{4} -\frac{k}{2N} \right\} =\frac{3}{4} -\frac{k}{2N}
>\frac{1}{2} > 1- 2\cd \frac{k}{N} ,$$
a contradiction.

If $N=14$, then $a=(3N/4)-(1/2)=10$ and $b=N-2-a=2$.
In this case $k=5$ is prime to $14$, $5/14<1/2$, but
$$\left\{ \frac{ka}{N}\right\} +\left\{ \frac{kb}{N}\right\} =
\left\{ \frac{50}{14}\right\} +\left\{ \frac{10}{14}\right\} =\frac{18}{14} >
1,$$
a contradiction.

If $N=18$, then $a=13$ and $b=3$. Then $k=5$ is prime to $18$, $5/18<1/2$, but
$$\left\{ \frac{ka}{N}\right\} +\left\{ \frac{kb}{N}\right\} =
\left\{ \frac{65}{18}\right\} +\left\{ \frac{15}{18}\right\}=\frac{26}{18} >1,$$
a contradiction.

If $N=22$, then $a=16$ and $b=4$. Then $k=5$ is prime to $22$, $5/22<1/2$, but
$$\left\{ \frac{ka}{N}\right\} +\left\{ \frac{kb}{N}\right\} =
\left\{ \frac{80}{22}\right\} +\left\{ \frac{20}{22}\right\}=\frac{34}{22} >
1,$$
a contradiction. Therefore, the case $a=(3N/4)-(1/2)$ is also impossible if
$N>10$. This completes the proof of the theorem. \hfill $\square$

Note that in the proof of Theorem \ref{t2} we only used Theorems
\ref{l1}-\ref{l3}. Therefore, as the proofs of Theorems \ref{l1}-\ref{l3} are
completed, Theorem \ref{t2} is also proved (subject to the number theoretic
Lemmas \ref{lL7}-\ref{lL1}).

\section{Proof of Theorem \ref{l4}}
By Theorem \ref{t2}, we may assume that the tiling is irregular.

By symmetry, we may assume $\be\le \ga$. Then, by $\al /\pi =(N-2)/(2N)$ we
have
$$\frac{\ga}{\pi}\ge \frac{\be +\ga}{2\pi}=\frac{\pi -\al}{2\pi}=\ha -
\left( \frac{1}{4} -\frac{1}{2N} \right) =\frac{1}{4}+
\frac{1}{2N} > \frac{1}{4}.$$
It follows that in every equation $p\al +q\be +
r\ga =v\pi$ we have $r\le 7$. Note that in every equation we have $p\le 4$,
as $\al >2\pi /5$ by $N>10$.

By the irregularity of the tiling,
there exists an equation $p_0 \al +q_0 \be +r_0\ga =v_0 \pi$ with $q_0 <r_0$.
We may assume $\min(p_0 ,q_0 )=0$, since otherwise we turn to the equation
$(p_0 -m) \al +(q_0 -m) \be +(r_0 -m)\ga =(v_0 -m) \pi$, where $m=\min(p_0 ,q_0 )$. We have
$$(p_0 -q_0 ) \al +(r_0 -q_0 )\ga =(v_0 -q_0) \pi .$$
We put
\begin{equation}\label{er46}
u=p_0 -q_0 , \quad  s=r_0 -q_0  , \quad t= 2v_0 -p_0 -q_0 .
\end{equation}
Note that $-6\le u\le 4$ and $1\le s\le 7$ by $p_0 \le 4$, $\min(p_0 ,q_0 )=0$
and $q_0 <r_0 \le 7$. It is clear that $t\le 4$.

By $u\al+s\ga =(v_0 -q_0) \pi$ we obtain
\begin{equation}\label{er44}
\ga =\frac{1}{s} \cd
\left[ v_0 -q_0 -u\cd \left(  \frac{1}{2} -\frac{1}{N} \right) \right] \pi =
\left[ \frac{t}{2s} +\frac{u}{sN}  \right] \pi
\end{equation}
and
\begin{equation}\label{er45}
\be =\pi -\al -\ga = 
\left[ \left( \frac{1}{2} -\frac{t}{2s} \right) +\frac{1}{N} -\frac{u}{sN} \right] \pi =
\left[ \frac{s-t}{2s} -\frac{u-s}{sN} \right] \pi .
\end{equation}
Since $\be >0$, we get $(s-t)N>2(u-s) =2(p_0 -r_0 )\ge -14$. Thus $s\ge t$, as
$s-t<0$ would imply $N<14$. Next we show $s\le 2t$. Suppose $s>2t$. Then
$$0\le (\ga -\be )/\pi =
\frac{2t-s}{2s} +\frac{2u-s}{sN} \le \frac{-1}{2s} +\frac{2u-s}{sN},$$
hence $1\le 2(2u-s)/N$, $N\le 2( 2u-s)\le 14$, which is impossible.
Thus $s\le 2t$, which also implies $t\ge 1$.

Summing up: we have
\begin{equation}\label{er47}
-6\le u\le 4, \quad 1\le s\le 7, \quad 1\le t\le 4 \quad \text{and} \quad t\le
s \le 2t.
\end{equation}
So the angles $\be$ and $\ga$ can only have a finite number (more precisely, at
most $11\cd 7\cd 4 =308$) of possible values for every $N$. We show that if $N
\ge 25$ and $N\ne 30, 42$, then only $\ga =\pi /2$ and $\ga =(1/2)-(1/N)$ are
possible, as the other cases do not satisfy Condition (K) and Condition (E).
We distinguish between two cases.

\msk \noi
Case I: $t=s$. By \eqref{er46}, this implies $2v_0 =p_0 +r_0$.
Then \eqref{er44} and \eqref{er45} give
$$\be = \frac{s-u}{sN} \cd \pi \quad \text{and} \quad
\ga =\left( \frac{1}{2} +\frac{u}{sN} \right) \cd \pi .$$
Then $\be >0$ gives $s>u$; that is, $r_0 >p_0$.

Thus the nonnegative integers $p_0 ,q_0 ,r_0 ,v_0$ satisfy the following
conditions: $v_0 =1$ or $2$, $\min(p_0 ,q_0 )=0$, $2v_0 =p_0 +r_0$, $p_0 <r_0$
and $q_0 <r_0$. It is easy to check that the quadruples $(p_0 ,q_0 ,r_0 ,v_0 )$
satisfying these conditions are the following:
\begin{align*}
&(0,0,2,1), \quad (0,1,2,1), \quad (0,0,4,2), \quad (0,1,4,2), \\
& (0,2,4,2), \quad (0,3,4,2)\quad  \text{and} \quad  (1,0,3,2).
\end{align*}
The values of $(s-u)/s =(r_0 -p_0 )/(r_0 -q_0 )$ obtained in these cases are
$1, 2, 4$, $2/3$ and $4/3$. That is, the possible values of $\be$ are $\pi /N$,
$2\pi /N$, $4\pi /N$, $2\pi /(3N)$ and $4\pi /(3N)$. The first two cases give
the triples listed in the theorem.

Suppose $\be =4\pi /N$. Then $\ga = ((1/2)-(3/N))\pi$. If $N\ge 43$, (ii) of
Lemma \ref{lL1} gives an integer $k$ such that $N/6< k<N/4$ and
$\gcd (k, 2N )=1$. Then $\{ kb/n\} =\{ 4k/N\} > 2/3$ and
$$\left\{ \frac{kc}{n} \right\} =
\left\{ \frac{k}{2} -\frac{3k}{N} \right\} = \left\{ \frac{1}{2} -
\frac{3k}{N} \right\} >\frac{3}{4} ,$$
since $1/2<3k/N<3/4$. Thus the triple $(\al ,\be ,\ga )$ does not satisfy
Condition (K). It is easy to check that for every $25\le N <42$
the triple $(\al ,\be ,\ga )=(((1/2)-(1/N))\pi , 4\pi /N, ((1/2)-(3/N))\pi )$
does not satisfy Condition (K).\footnote{In this computation and also in
the computer search needed in the proof of the next theorem I applied GNU Octave (https://www.gnu.org/software/octave/).}

Therefore, the case $\be =4\pi /N$ is impossible if $N\ge 25$
and $N\ne 42$.

Next suppose $\be =2\pi /(3N)$. Then $\ga =(1/2)+(1/(3N))$. Let
$$ k=
\begin{cases}
N+1 &\text{if $N\equiv 0$ or $4$ (mod $6$)},\\
N+2 &\text{if $N\equiv 3$ or $5$ (mod $6$)},\\
N+3 &\text{if $N\equiv 2$ (mod $6$)},\\
N+4 &\text{if $N\equiv 1$ (mod $6$)}.
\end{cases}
$$
Then $\gcd (k,6N)=1$, and $\{ k/N\} <1/2$. We have $\{ kb/n \} = \{ 2k/(3N) \} >
2/3$ and
$$\left\{ \frac{kc}{n} \right\} =\left\{ \frac{k}{2} +\frac{k}{3N} \right\}
=  \left\{ \frac{1}{2} +\frac{k}{3N} \right\} >\frac{5}{6},$$
since $1/3 <\{ k/(3N)\} <1/2$. Thus $(\al ,\be ,\ga )$ does not satisfy
Condition (K), and the case $\be = 2\pi /(3N)$ is impossible. 

Finally, suppose $\be =4\pi /(3N)$. Then $\ga = (1/2)-(1/(3N))$. We put
$k=2N+1$ if $N \not\equiv 1$ (mod $3$), and $k=2N+3$ if $N \equiv 1$ (mod $3$).
Then $\gcd (k,6N)=1$ and $\{ k/N\} <1/2$.
We have $\{ kb/n\} =\{ 4k/(3N)\} \ge 2/3$, since $8/3 < 4k/(3N)<3$. On the
other hand,
\begin{equation}\label{er42}
\left\{ \frac{kc}{n} \right\} =\left\{ \frac{k}{2} - \frac{k}{3N} \right\} 
> \frac{3}{4} ,
\end{equation}
since $2/3< k/(3N)<3/4$. Thus $(\al ,\be ,\ga )$ does not satisfy Condition
(K). Therefore, the case $\be =4\pi /(3N)$ is also impossible.

\msk \noi
Case II: $t<s$. First suppose $N>500$. Then, by \eqref{er45} we have 
$$\frac{\be}{\pi} = \frac{s-t}{2s} -\frac{u-s}{sN}  \ge \frac{1}{2s} -
\frac{3}{sN} =\frac{N-6}{2sN} \ge \frac{N-6}{14N} >\frac{1}{15}.$$
This implies that $q<30$ holds in every equation $p\al +q\be +r\ga =v\pi$.

Let $p\al +q\be +r\ga =v\pi$ be any of these equations.
Substituting \eqref{er44} and \eqref{er45} into $p\al +q\be +r\ga =v\pi$ we
obtain
$$p\cd \ha +q \left( \frac{1}{2} -\frac{t}{2s} \right) + r\cd \frac{t}{2s}
+ \frac{1}{N} \cd \left[ -p +q \left( 1 -\frac{u}{s} \right) + r\cd \frac{u}{s}
  \right] =v$$
and $A\cd N=2\cd (-ps+q(s-u)+ru)$, where $A=2sv-(ps+q(s-t)+rt)$. If $A\ne 0$,
then
$$N\le 2\cd |-ps+q(s-u)+ru|\le 2\cd \max(qs+ru, ps+qu)\le 2\cd (30\cd 7 +7\cd 
4)<500,$$
which is impossible. Therefore, we have $A=0$, hence
$-ps+q(s-u)+ru =0$.

We proved that $-ps+q(s-u)+ru =0$ holds for every equation $p\al +q\be +r\ga
=v\pi$. Let $K$ denote the number of the tiles. Taking the sum of the equations
$-ps+q(s-u)+ru =0$ we obtain $0=-(K -2N)s +K(s-u) +Ku =2NS$, a contradiction.
Therefore, Case II is impossible if $N>500$.

If $N\le 500$, then we check for every possible triple $(\al ,\be ,\ga )$
whether or not it satisfies Condition (K) and Condition (E). If $N$ is given,
then $\be$ and $\ga$ are determined by \eqref{er45} and \eqref{er44}. As these
formulas show, we may take $n=2sN$. We check, for every choice of $u,s,t$
satisfying \eqref{er47} and also $t<s$ whether or not \eqref{e13} holds for
every $k$ such that $\gcd (k,nN)=1$ and $\{ k/N\} <1/2$.

A computer search shows that in the range $60 <N\le 500$ only $N=78$ produces
triples $(\al ,\be ,\ga )$ satisfying Condition (K).
More precisely, for $N=78$ there is just one such triple, namely
$$\left( \frac{38}{78} \pi , \frac{17}{78} \pi , \frac{23}{78} \pi \right) .$$
However, the only equations $p\al +q\be +r\ga =v\pi$ in this case
are $\al +\be +\ga =\pi$ and $2\al +2\be +2\ga =2\pi$. Thus (iii) of Condition
(E) is not satisfied, since we have $p>q$ in the equations at the vertices
$V_j$. Thus the case $N=78$ cannot occur.

In the range $42 <N\le 60$ only $N=60$ produces
triples $(\al ,\be ,\ga )$ satisfying Condition (K).
For $N=60$ there are two such triples, namely
\begin{equation}\label{er50}
\left( \frac{29}{60} \pi , \frac{12}{60} \pi , \frac{19}{60} \pi \right)
\quad \text{and} \quad
\left( \frac{29}{60} \pi , \frac{11}{60} \pi , \frac{20}{60} \pi \right) .
\end{equation}
In the first case the only equations $p\al +q\be +r\ga =v\pi$ 
are $5\be =\pi$, $10\be =2\pi$, $\al +\be +\ga =\pi$, $2\al +2\be +2\ga =2\pi$
and $\al +6\be +\ga =2\pi$.

We can see that $q\ge r$ holds in each of these equations. Then (iii) of
Condition (E) can hold only if the equations with $q>r$ do not occur in the
tiling. The remaining equations are $\al +\be +\ga =\pi$ and
$2\al +2\be +2\ga =2\pi$. Thus Condition (E) is not
satisfied, since we have $p>q$ in the equations at the vertices
$V_j$. Thus this case is impossible. 

If $(\al ,\be ,\ga )$ equals the second triple of \eqref{er50}, then 
the equations $p\al +q\be +r\ga =v\pi$ 
are the following: $3\ga =\pi$, $6\ga =2\pi$, $\al +\be +\ga =\pi$, $2\al +2\be +2\ga =2\pi$, $\al +\be +4\ga =2\pi$,
and $3\al +3\be  =2\pi$.

We can see that $p=q$ holds in each of these equations. Since
$p>q$ holds in the equations at the vertices $V_j$, Condition (E) is not
satisfied, and this case is also impossible. 

In the range $24<N\le 42$ only $N=30$ and $N=42$ produce
triples $(\al ,\be ,\ga )$ satisfying Condition (K). This completes the
proof of the theorem. \hfill $\square$

\section{Proof of Lemma \ref{lL7}}
We may assume $\gcd (N',2nN)=1$, since otherwise we replace $N'$ by $N'+jN$ with
a suitable $j$.
 
Suppose there is an odd prime $p$ such that $p\mid n$
and $p\ {\mid \!\!\! \not} \quad N$. Let $P$ denote the product of primes
dividing $n$ and different from $p$. (Put $P=1$ if there is no such prime.)
Let $(NPa)/n=M/m$, where $\gcd (M,m)=1$.
Since $p\ {\mid \!\!\! \not} \quad NPa$ and $p\mid n$, we have $p\mid m$,
and thus $m\ge p\ge 3$.

Let $s$ be such that $s M\equiv 1$ (mod $m$). Then $p\ {\mid \!\!\! \not} \quad
s$, as $p\mid m$. Put $k_i =N' +isNP$ for every integer $i$. Then $k_i$ is
not divisible by any prime divisor of $nN$ except perhaps $p$. But if $p\mid
k_i$, then $p\ {\mid \!\!\! \not} \quad k_{i-1} ,k_{i+1}$, since
$p\ {\mid \!\!\! \not} \quad sNP$. Thus either $k_{i}$ is prime to $nN$ or
both of $k_{i-1} ,k_{i+1}$ are prime to $nN$. Now
$$\frac{k_i a}{n} = \frac{N'a}{n} +i\frac{sNPa}{n} =
\frac{N'a}{n} +i\frac{sM}{m} \equiv \frac{N'a}{n} +\frac{i}{m}
\qquad \text{(mod $1$)}.$$
This implies, by $m\ge 3$, that there are two consecutive $i$'s with
$\{ k_i a/n\} \ge 1/3$. For at least one of them, $k_i$ is prime to $nN$.
We find that (i) holds.

Next suppose that every odd prime divisor of $n$ divides $N$.
Suppose $N$ is odd. Then $k_i =N' +2iN$ is prime to $nN$ for every $i$. Now
$k_i a/n = (N'a/n) +i\cd (2N/n)$ and thus, if $n\ {\mid \!\!\! \not} \quad 2N$,
then for a suitable $i$ we have $\{ k_i a/n\} \ge 1/2$. That is, we have either
(i) or (ii) in this case.

If $N$ is even, then $k_i =N' +iN$ is prime to $nN$ for every $i$. Since
$k_i a/n = (N'a/n) +i\cd (N/n)$, we find that if
$n\ {\mid \!\!\! \not} \quad N$, then for a suitable $i$
we have $\{ k_i a/n\} \ge 1/2$. That is, we have either (i) or (iii) in this
case. This completes the proof. \hfill $\square$

\section{Proof of Lemma \ref{lL11}}

By symmetry, we may assume $p\ge q$. Let $a/n=a_1 /n_1$ and $b/n=b_2 /n_2$,
where $\gcd (a_1 ,n_1 )=\gcd (b_2 ,n_2 ) =1$. Applying \eqref{eL11} with $k=1$
we obtain
\begin{equation}\label{eL13}
\frac{pa+qb}{n} =\frac{N-2}{N} . 
\end{equation}
We consider three cases.

\noi
Case I: $N$ is odd. Then $N'=(N-1)/2$ is prime to $N$. Suppose $n_1
\mid 2N$. For a suitable $j$, $k_1 =((N-1)/2) +jN$ is prime to $nN$. By
\eqref{eL11} we obtain
$$ \frac{p}{2N} +\ep \le \frac{p}{n_1} +\ep  \le p\cd \left\{ \frac{k_1 a_1 }{n_1} \right\} +\ep = p\cd \left\{ \frac{k_1 a}{n} \right\} +\ep = 1-2\left\{ \frac{k_1}{N} \right\}   =\frac{1}{N},$$
where $\ep =q\cd \{ k_1 b/n\}$. Therefore, we have $p\le 2$. If $p=2$, then
$\ep =0$ and $q=0$. If $p\le 1$, then $q\le 1$, and we have $p+q\le 2$ in both
cases.

Therefore, we may assume that $n_1$ does not divide $2N$.
Then, applying Lemma \ref{lL7}, we find
that (i) of Lemma \ref{lL7} holds with $a_1$ and $n_1$ in place of $a$ and
$n$. That is, there is a $k$ prime to $n_1 N$ such that $k\equiv (N-1)/2$
(mod $N$) and $\{ ka_1 /n_1 \}\ge 1/3$. For a suitable $j$,
$k_2 =k+jn_1 N$ will be prime to $nN$. Then \eqref{eL11} gives
$$ \frac{p}{3} \le p\cd \left\{ \frac{k_2 a_1 }{n_1} \right\} = 
p\cd \left\{ \frac{k_2 a}{n} \right\} \le 1-2\left\{ \frac{k_2}{N} \right\} = \frac{1}{N} \le \frac{1}{3},$$
$p\le 1$, and we are done.

\noi
Case II: $4\mid N$. Then $N'=(N/2)-1$ is prime to $N$. Suppose $n_1
\mid N$. For a suitable $j$, $k_3 =(N/2)-1 +jN$ is prime to $nN$. By
\eqref{eL11} we obtain
$$ \frac{p}{N} +\ep  \le \frac{p}{n_1} +\ep \le p\cd \left\{ \frac{k_3 a_1 }{n_1} \right\} +\ep = 1-2\left\{ \frac{k_3}{N} \right\}   =\frac{2}{N},$$
where $\ep =q\cd \{ k_3 b/n\}$. From this we obtain $p+q\le 2$ as in case I.

If $n_1 \ {\mid \!\!\! \not} \quad N$, then applying Lemma
\ref{lL7}, we find that (i) of Lemma \ref{lL7} holds
with $a_1$ and $n_1$ in place of $a$ and $n$. That is, there is
a $k$ prime to $n_1 N$ such that $k\equiv (N/2)-1$ (mod $N$) and
$\{ ka_1 /n_1 \}\ge 1/3$. For a suitable $j$,
$k_4 =k+jn_1 N$ will be prime to $nN$. Then \eqref{eL11} gives
$$ \frac{p}{N}+\ep \le \frac{p}{3}+\ep \le p\cd \left\{ \frac{k_4 a_1 }{n_1} \right\} +\ep = 1-2\left\{ \frac{k_4}{N} \right\} = \frac{2}{N},$$
where $\ep =q\cd \{ k_4 b/n\}$. From this inequality we obtain $p+q\le 2$
as above.

\noi
Case III: $N$ is even and $N/2$ is odd. Then $N'=(N/2)-2$ is prime to $N$.
Note that the first possible value of $N$ is $10$, as $N=6$ is excluded.

\noi
Case IIIa: $n_1 \mid N$. Then $a/n=u/N$, where $0<u<N$ is an 
integer. Suppose $p\ge 2$. For a suitable $j$, $k_5 =(N/2)-2 +jN$
is prime to $nN$. By \eqref{eL11} we obtain
\begin{equation}\label{e3}
p\cd \left\{ \frac{k_5 u }{N} \right\} + q\cd
\left\{ \frac{k_5 b}{n} \right\} = 1-2\left\{ \frac{k_5}{N} \right\}
=\frac{4}{N}.
\end{equation}
Since $\{ k_5 u/N\}$ is a positive integer multiple
of $1/N$ and $p\ge 2$, we have $\{ k_5 u/N\} =1/N$ or $2/N$.
If $\{ k_5 u/N\} =2/N$, then \eqref{e3} gives $p=2$, $q=0$, and we are done.

If $\{ k_5 u/N\} =1/N$, then
$k_5 u\equiv 1$ (mod $N$), $u$ is odd, $(N/2)u \equiv N/2$ (mod $N$),
$k_5 u\equiv ((N/2)-2)u \equiv (N/2)-2u \equiv 1$ (mod $N$) and
$u\equiv ((N/2)-1)/2$ (mod $N/2$).
Now $2u/N=2a/n<1$ by \eqref{eL13}, and thus $u=((N/2)-1)/2=(N-2)/4$.

Since $N\ge 10$, $(N/2)-4$ is also prime to $N$. For a suitable $j$,
$k_6 = (N/2)-4+jN$ is prime to $nN$. Then we have
$$p\left\{ \frac{k_6 u}{N} \right\} + q\left\{ \frac{k_6 b}{n} \right\} =
1-2\left\{ \frac{k_6}{N} \right\}   =\frac{8}{N},$$
and thus $\{ k_6 u/N\} \le 4/N$. However, we have
$$k_6 u\equiv k_5 u -2u\equiv 1-2u= 1-\frac{N-2}{2} \equiv \frac{N+4}{2}
\qquad \text{(mod $N$)}$$
and $\{ k_6 u/N\} =(N+4)/(2N)>1/2>4/N$, a contradiction. Therefore, we have
$p\le 1$ and $p+q\le 2$.

\noi
Case IIIb: $n_1 \ {\mid \!\!\! \not} \quad N$. By Lemma \ref{lL7}, there is a
$k$ prime to $n_1 N$ such that $k\equiv (N/2)-2$ (mod $N$) and
$\{ ka_1 /n_1 \}\ge 1/3$. For a suitable $j$, $k_8 =k+jn_1 N$ will be prime to
$nN$. Then \eqref{eL11} gives
$$ \frac{p}{3} \le p\cd \left\{ \frac{k_8 a_1 }{n_1}
\right\} +q\cd \left\{ \frac{k_8 b_2 }{n_2} \right\}   =
1-2\left\{ \frac{k_8}{N} \right\} = \frac{4}{N} \le \frac{4}{10} <
\frac{2}{3}.$$
Thus $p\le 1$, $p+q\le 2$, and the proof is complete. \hfill $\square$

\section{Proof of Lemma \ref{lL1}}

\begin{lemma}\label{lL2}
Let $u,m,N$ be integers such that $m,N>0$ and $\gcd (u,m)=1$. Let
$p_1 \stb p_s$ be those primes that divide $N$ but not $m$. If $c>0$ and
\begin{equation}\label{eL14}
\frac{cN}{m} \left( 1-\frac{1}{p_1} \right) \ldots  \left( 1-\frac{1}{p_s}
\right) \ge 2^s ,
\end{equation}
then for every real number $a$ there is an integer $k$ such that $a\le k<a+cN$,
$k\equiv u$ (mod $m$), and $\gcd (k,N)=1$.
\end{lemma}
\proof
Let $A_d$ denote the set of integers $k$ such that $a\le k<a+cN$, $k \equiv
u$ (mod $m$) and $d\mid k$. If $\gcd (d,m)=1$, then there is a $j_0$ such that
$j_0 m\equiv -u$ (mod $d$), and then $A_d$ equals the set of numbers
$u+j_0 m+ jmd$ such that $a\le u+j_0 m+jmd <a+cN$. Thus $|A_d |$ equals the
number of integers $j$ with $b\le j<b+(cN/md)$, where $b=(a-u-j_0 m)/(md)$. 
Therefore, we have $|A_d |=(cN/md)+\ep _d$, where $|\ep _d |<1$. 

If $S$ denotes the number of integers $k$ such that $a\le k<a+cN$,
$k\equiv u$ (mod $m$), and $\gcd (k,N)=1$, then
\begin{align*}
S &=\sum_{d\mid p_1 \cdots p_s} \mu (d) |A_d | =\sum_{d\mid p_1 \cdots p_s} \mu (d)
\frac{cN}{md} +\sum_{d\mid p_1 \cdots p_s} \mu (d)\cd \ep _d >\\
&> \frac{cN}{m} \left( 1-\frac{1}{p_1} \right) \ldots  \left( 1-\frac{1}{p_s}
\right) - 2^s .
\end{align*}
If \eqref{eL14} is true, then $S>0$, which proves the lemma. \hfill $\square$

\bsk
{\bf Proof of (i) of Lemma \ref{lL1}.} Let $p_1 \stb p_s$ be the odd prime
divisors of the even number $N$. By Lemma \ref{lL2}, statement (i) of Lemma
\ref{lL1} is true, if
\begin{equation*}
\frac{N}{16} \left( 1-\frac{1}{p_1} \right) \cdots  \left( 1-\frac{1}{p_s}
\right) \ge 2^s . 
\end{equation*}
If $s\ge 4$, then
$$N\cd \( 1-\frac{1}{p_1} \) \cdots  \( 1-\frac{1}{p_s} \) \ge 2\cd (p_1 -1)
\cdots (p_s -1)\ge 2\cd 2\cd 4\cd 6 \cd 10^{s-3} >16\cd 2^s ,$$
and thus the statement is true. Therefore, we may assume $s\le 3$. If
$N>480$, then
$$\frac{N}{16} \left( 1-\frac{1}{p_1} \right) \cdots  \left( 1-\frac{1}{p_s}
\right) > \frac{480}{16} \cd \frac{1}{2} \cd \frac{2}{3} \cd \frac{4}{5} =
2^3 ,$$
and then the statement is true again. Finally, it is easy to check that for
every even integer $N\in [26 ,480]$ there are integers $k, k'$ with the required
properties.

\bsk
{\bf Proof of (ii) of Lemma \ref{lL1}.} Let $p_1 \stb p_s$ be the odd prime
divisors of $N$. Applying Lemma \ref{lL2} with $m=2$ and $u=1$ we obtain that
statement (ii) of Lemma \ref{lL1} is true, if
\begin{equation*}
\frac{N}{24} \left( 1-\frac{1}{p_1} \right) \cdots  \left( 1-\frac{1}{p_s}
\right) \ge 2^s . 
\end{equation*}
If $s\ge 4$, then
$$N\cd \( 1-\frac{1}{p_1} \) \cdots  \( 1-\frac{1}{p_s} \) \ge (p_1 -1)\cdots
(p_s -1)\ge 2\cd 4\cd 6 \cd 10^{s-3} >24\cd 2^s ,$$
and thus the statement is is true. Therefore, we may assume $s\le 3$. If
$N>720$, then
$$\frac{N}{24} \left( 1-\frac{1}{p_1} \right) \cdots  \left( 1-\frac{1}{p_s}
\right) > \frac{720}{24} \cd \frac{1}{2} \cd \frac{2}{3} \cd \frac{4}{5} =
2^3 ,$$
and then the statement is true again. Finally, it is easy to check that for
every integer $N\in [43 ,720]$ there is an integer $k$ with the required
properties. \hfill $\square$

\hfill \eject

\noindent
\textsc{E\"otv\"os Lor\'and University (ELTE), Department of Analysis}
 
\noindent
\textsc{P\'az\-m\'any P\'e\-ter s\'et\'any 1/c, H-1117, Budapest, Hungary}

\noindent
\textit{Email address}: {\tt miklos.laczkovich@gmail.com}


\begin{thebibliography}{100}

\bibitem{F} G. Frederickson: {\it Dissections: Plane and Fancy.} Cambridge
Univ. Press, 1997.

\bibitem{G} B. Gr\"unbaum and C. C. Shephard: {\it Tilings and Patterns.}
W. H. Freeman, New York, 1987. 

\bibitem{L} M. Laczkovich, Tilings of polygons with similar triangles,
\it Combinatorica \rm 10 (1990), 281-306.

\bibitem{L2} M. Laczkovich, Tilings of polygons with similar triangles II,
{\it Discrete and Computational Geometry,} {\bf 19} (1998), 411-425.

\bibitem{L3} M. Laczkovich, Tilings of convex polygons with congruent triangles,
{\it Discrete and Computational Geometry,} {\bf 48} (2) (2012), 330-372.

\bibitem{L4} M. Laczkovich, Some tilings of the regular pentagon and octagon,
in preparation.

\bibitem{La} C. D. Langford, Tiling patterns for regular polygons,
{\it The Mathematical Gazette} {\bf 44} (1960), 105-110. 
  
\bibitem{Sz} B. Szegedy, Tilings of the square with similar right triangles,
{\it Combinatorica} {\bf 21} (2001), no. 1, 139-144. 


\end{thebibliography}
\end{document}